\documentclass[reqno,12pt]{amsart}
\usepackage{amsmath,amsthm,amscd,amssymb}
\newcommand{\bbR}{{\mathbb{R}}}

\newcommand{\calM}{{\mathcal M}}

\newcommand{\lb}{\label}
\newcommand{\f}{\frac}

\newcommand{\loc}{\text{\rm{loc}}}
\newcommand{\spec}{\text{\rm{spec}}}

\newcommand{\bi}{\bibitem}

\newcommand{\beq}{\begin{equation}}
\newcommand{\eeq}{\end{equation}}
\newcommand{\ba}{\begin{align}}
\newcommand{\ea}{\end{align}}

\newcounter{smalllist}
\newenvironment{SL}{\begin{list}{{\rm\roman{smalllist})}}{%
\setlength{\topsep}{0mm}\setlength{\parsep}{0mm}\setlength{\itemsep}{0mm}%
\setlength{\labelwidth}{2em}\setlength{\leftmargin}{2em}\usecounter{smalllist}%
}}{\end{list}}

\allowdisplaybreaks

\newtheorem{theorem}{Theorem}

\theoremstyle{definition}

\theoremstyle{remark}

\newcommand{\abs}[1]{\lvert#1\rvert}

\begin{document}

\title[Connectedness of the Isospectral Manifold]
{Connectedness of the Isospectral Manifold for One-Dimensional Half-Line
Schr\"odinger Operators}
\author[F. Gesztesy and B. Simon]{Fritz Gesztesy$^{1}$ and Barry Simon$^{2}$}

\thanks{$^1$ Department of Mathematics, University of Missouri,
Columbia, MO 65211, USA.
E-mail: fritz@math.missouri.edu}
\thanks{$^2$ Mathematics 253-37, California Institute of Technology,
Pasadena, CA 91125,
USA. E-mail: bsimon@caltech.edu. Supported in part by NSF grant DMS-0140592}


\dedicatory{To Elliott Lieb on his 70th birthday, with our best wishes}

\date{June 19, 2003}

\begin{abstract} Let $V_0$ be a real-valued function on $[0,\infty)$
and $V\in L^1([0,R])$ for all $R>0$ so that $H(V_0)=  -\f{d^2}{dx^2}+V_0$
in $L^2([0,\infty))$ with $u(0)=0$ boundary conditions has discrete
spectrum bounded from  below. Let $\calM (V_0)$ be the set of $V$ so that
$H(V)$ and $H(V_0)$ have the same spectrum.  We prove that $\calM(V_0)$
is connected.
\end{abstract}

\maketitle

The limitation of our knowledge of inverse spectral theory for
Schr\"odinger operators
$H_\bbR(V)= -\f{d^2}{dx^2} + V$ in $L^2(\bbR)$ is shown by the
following open question:  Let $\calM_\bbR (x^2)$ denote the set of all
$V$'s with $\sigma(H_\bbR (V))=\{1,3,5,\dots\}$,  the spectrum of the
harmonic oscillator $-\f{d^2}{dx^2}+x^2$ in $L^2(\bbR)$. Is $\calM_\bbR
(x^2)$ connected?  One can ask the same question for restricted sets of
$V$'s, say requiring that $V$ is $C^k(\mathbb{R})$ for some
$k\in\mathbb{N}$ or $C^\infty(\mathbb{R})$.  The question remains open for
sets so large that they include the translates of $V$\!.

If one demands that $V$ be close to $x^2$ in a strong sense, there
are some results that go
back to McKean--Trubowitz \cite{MT} (see also Levitan \cite{Lev}),
culminating in the recent
beautiful paper of Chelkak, Kargaev, and Korotyaev \cite{CKK} who
require $V(x)=x^2 + q(x)$
with $\int [\abs{q'(x)}^2 + x^2 \abs{q^2(x)}^2]\, dx <\infty$ (by a
Sobolev estimate such
$q$'s have $\abs{q(x)}\to 0$ as $\abs{x}\to\infty$). Their analysis
implies the set of $V\in
\calM_\bbR (x^2)$ obeying this estimate is connected.

The purpose of this note is to make what turns out to be an elementary
observation: The corresponding problem for the half-line is easy! Suppose
$V$ is real-valued and in $L^1([0,R])$ for all $R>0$. We consider
$H(V)=-\f{d^2}{dx^2} + V$ in $L^2  ([0,\infty))$ with
$u(0)=0$ (i.e., Dirichlet) boundary conditions. If $V$ is merely in
$L^1([0,R])$ for all $R>0$, we will use that to define a topology
on $\calM(V)$. If $V$ is in $C^k([0,\infty))$ for some $k\in\mathbb{N}$,
we will use the norm on $C^k([0,R])$ for all $R>0$ to define the
topology on $\calM(V)$. Here is our result:

\begin{theorem} \lb{T1} Let $V_0, V_1$ be real-valued and in $L^1([0,R])$
for all $R>0$ so that $H(V_\ell)$ is bounded  from below {\rm{(}}which
means that $-\f{d^2}{dx^2} +V_\ell$ is in the limit point case at 
$\infty$, cf.\
\cite{Har}{\rm{)}} and each has discrete spectrum $\spec
(H(V_\ell))=\{E_j(V_\ell)\}_{j\in\mathbb{N}}$, $\ell=0,1$.

Suppose
\begin{equation} \lb{1.1}
\spec(H(V_0)) = \spec(H(V_1)).
\end{equation}
Then there exists $\{V_t\}_{0\leq t\leq 1}$, with $V_t$ real-valued and in
$L^1([0,R])$ for all $R>0$, interpolating $V_0$ and $V_1$ so that
\[
\spec(H(V_t))=\spec(H(V_0)), \quad t\in [0,1]
\]
and $t\to V_t|_{[0,R]}$ is continuous in $L^1([0,R])$ for all $R>0$.
Moreover,
\begin{SL}
\item[{\rm{(i)}}] $t\to V_t|_{[0,R]}$ is real analytic in $L^1([0,R])$
for all $R>0$.
\item[{\rm{(ii)}}] If $V_0$ and $V_1$ are $C^k([0,\infty))$ for some
$k\in\mathbb{N}$, then $V_t$ is
$C^k([0,\infty))$ and $t\to V_t|_{[0,R]}$ is real analytic in
$C^k([0,R]))$ for all $R>0$.
\end{SL}
\end{theorem}

We recall that all eigenvalues $\{E_j(V)\}_{j\in\mathbb{N}}$ of the
Dirichlet operator $H(V)$ in $L^2([0,\infty))$ are simple.

The proof exploits the $A$-function studied by us in \cite{GS,S}. $A$
can be defined in terms
of the spectral measure defined in the standard way (see, e.g.,
\cite{Mar}) so that the Weyl
$m$-function satisfies
\begin{equation} \lb{1.2}
m(z) = c + \int_\bbR  \left[ \f{1}{\lambda-z} -\f{\lambda}{1+\lambda^2}
\right] d\rho(\lambda), \quad z\in\mathbb{C}\backslash\spec(H(V)).
\end{equation}
We let $d\rho_0$ denote the spectral measure for the case $V=0$. It
is well-known, since
\begin{equation} \lb{1.3}
m(z,V=0)=i\sqrt{z},
\end{equation}
that
\begin{equation} \lb{1.4}
d\rho_0 (E) = \f{1}{\pi} \sqrt{E}\, \chi_{[0,\infty)}(E)\, dE.
\end{equation}
$A$ is then defined by
\begin{equation} \lb{1.1a}
A(\alpha) = -2 \int_{-\infty}^\infty \lambda^{-\f12} \sin(2\alpha
\sqrt{\lambda}\,) \, [d\rho(\lambda)-d\rho_0(\lambda)],
\end{equation}
where the integral is intended in the distributional sense on
$(-\infty,\infty)$ (so, a priori, $A$  is only a distribution, not a
function). Of course, $A = 0$ for $V=0$.

One can also define a distribution $A$ by \eqref{1.1a} with $\rho_0$
dropped (cf.\ \cite{GS}) but then this distribution is only $L_\loc^1$
away from $0$ and has a $\delta'$ singularity at $0$. We will not use
this approach in this note.

The fact on which Theorem~1 depends is the following:

\begin{theorem}\lb{T2} $d\rho$ is the spectral measure of an $H(V)$ with
$V\in L^1([0,R])$ for all $R>0$, $V$ real-valued, if and only
if $A\in L_\loc^1 (\mathbb{R})$. $V\in C^k([0,\infty))$ for some
$k\in\mathbb{N}$ if and only if $A\in C^k (\mathbb{R})$. If $d\rho_t$ is a
family so that $A_t|_{[-R,R]}$ is real analytic in $t$ in $L^1([-R,R])$
{\rm{(}}resp.~$C^k([-R,R])${\rm{)}} for all $R>0$, then $t\to V_t$ is real
analytic in $t$ in $L^1([0,R])$ {\rm{(}}resp.~$C^k([0,R])${\rm{)}} for all
$R>0$.
\end{theorem}

This theorem combines results from Gesztesy--Simon \cite{GS} and Simon
\cite{S} (who show $V$  is $C^k([0,\infty))$ if and only if $A$ is, once
one knows $V$ exists) and work of Remling \cite{Rem} or a suitable version
of the  Gel'fand--Levitan theory \cite[Ch.\ 2]{Lev87}, \cite[Sect.\
2.3]{Mar} to get the existence part of Theorem~\ref{T2}.

Once one has Theorem~\ref{T2}, Theorem~\ref{T1} is immediate.

\begin{proof}[Proof of Theorem~\ref{T1}] Let $\{E_j\}_{j\in\mathbb{N}}$ be
the common spectrum of $V_0$ and $V_1$ so
\begin{equation} \lb{1.5}
d\rho_\ell (E)= \sum_{j\in\mathbb{N}} a_{j,\ell} \, \delta (E-E_j), \quad
\ell=1,2.
\end{equation}
Define
\[
d\rho_t (E)= \sum_{j\in\mathbb{N}} \, [ta_{j,1} + (1-t) a_{j,0}]\delta
(E-E_j), \quad t\in [0,1].
\]
The associated $A$-functions satisfy
\begin{equation} \lb{1.6}
A_t = tA_1 + (1-t)A_0, \quad t\in [0,1].
\end{equation}
Clearly, if $A_0, A_1$ are $L_\loc^1(\mathbb{R})$, so is $A_t$, and if
$A_0, A_1$ are $C^k (\mathbb{R})$, so are $A_t$. Thus,
$d\rho_t$ is the spectral measure of a potential $V_t$ that, by
Theorem~\ref{T2}, has the claimed  properties.
\end{proof}

{\it Remarks.} 1. The $C^k$ result extends to $C^\infty$ and one can also
extend it to real analyticity.

\smallskip
2. Because $V$ on $[0,x_0]$ only depends on $A$ on $[0,x_0]$ (see
\cite{GS,S}), if $V_0$ and $V_1$
are $C^k$ on $[0,x_0]$, so is each $V_t$.

\smallskip
The key to our proof is the fact that, while not all measures
$\sum_{j\in\mathbb{N}} a_j\, \delta (E-E_j)$ are
spectral measures (the fact that the $A$-transform of $\rho-\rho_0$ has
no singularity at $\alpha=0$, which means the $\rho$-term alone  has a
specific singularity that places restrictions  on the $a_j$), those that
are form a convex, hence, connected set.

The difficulty of extending this to potentials on $\bbR$ is that there is
no known way to describe when a candidate spectral measure is, in fact,
the spectral  measure for a potential on $\bbR$. Typically, one reduces
the whole-line problem on $\mathbb{R}$ to two half-line problems on
$(-\infty,x_0)$ and $(x_0,\infty)$ for some fixed $x_0\in\bbR$ (cf.\
\cite[Ch.\ 7]{Lev87}, \cite{Rof67, Rof91}) but, in general, loses
control over the potential at the point $x_0$ (in the sense that
generally the potential $V$ will be discontinuous at $x=x_0$).  To
determine the potential, the spectral measure is a $2\times 2$ matrix in
this case \cite[Sect.\ 9.5]{CL}, \cite[Sect.\ X.III.5]{DS},  and
because of restrictions on this matrix, convex combinations of the matrix
measures coming  from a potential will {\it not} come from a potential, so
our method cannot extend.

We, of course, believe that for the whole-line case, there is also a
result on connectedness of
the spectral manifold for a potential like $V(x)=x^2$. But we wonder
about a case like
\begin{equation} \lb{1.7}
V(x) = \ln (\abs{x}+1) + \exp(x)
\end{equation}
with very different asymptotics if $x\to\infty$ and $x\to -\infty$.
$V(x)$ and $V(-x)$ are
obviously isospectral, but we wonder if there is a path between them
in the isospectral manifold.
It might be that the correct conjecture is that for any $V$\!,
$\calM_\bbR(V)$ consists of either one or two connected components.

We hope that this note will stimulate more work on this problem.

\smallskip
We close by thanking Elliott Lieb for his many years of research
results, insights, and service
to the mathematics and physics communities, and by expressing the
hope that he has enjoyed this
birthday bouquet.

\bigskip


\end{document}